\documentclass{IEEEtran}
\usepackage{cite}
\usepackage{amsmath,amssymb,amsfonts}
\usepackage{algorithmic}
\usepackage{graphicx}
\usepackage{textcomp}

\usepackage{amsthm}
\def\BibTeX{{\rm B\kern-.05em{\sc i\kern-.025em b}\kern-.08em
    T\kern-.1667em\lower.7ex\hbox{E}\kern-.125emX}}

\newtheorem{theorem}{Theorem}[section]
\newtheorem{corollary}[theorem]{Corollary}

\newtheorem{remark}[theorem]{Remark}
\newtheorem{assumption}[theorem]{Assumption}

\newcommand{\SO}{\operatorname{SO}(3)}
\newcommand{\SE}{\operatorname{SE}(3)}

\newcommand{\Lalge}{{\mathfrak g}}
\newcommand{\GLn}{{\operatorname{GL}(n)}}
\newcommand{\proj}{\operatorname{\pi_{\Lalge}}}
\newcommand{\meas} {{F}}
    
\begin{document}
\title{Design of Globally Exponentially Convergent Continuous Observers for Velocity Bias and State for    Systems on Real Matrix  Groups}
\author{Dong Eui Chang
\thanks{}}
\date{August 7, 2018}
\maketitle

\begin{abstract}
We propose globally exponentially convergent continuous observers for   invariant kinematic  systems on finite-dimensional matrix Lie groups. Such an observer  
estimates, from measurements of landmarks, vectors and biased velocity,  both the system state and   the unknown constant bias in velocity measurement, where the state belongs to the state-space Lie group  and the velocity  to the Lie algebra of the Lie group.  The main technique is to embed a given system defined on a matrix Lie group into  Euclidean space and build observers in the Euclidean space. The theory is illustrated with the special Euclidean group in three dimensions.
\end{abstract}


\section{Introduction} 
Consider an invariant kinematic  system on a matrix Lie group $G$:
\[
\dot g = g \xi
\]
with  $(g,\xi) \in G \times \mathfrak{g}$, where $G$ is embedded in $\mathbb R^{n\times n}$ and  $\mathfrak{g}$ denotes the Lie algebra of $G$. Suppose that  the velocity $\xi$ is measured with an additive unknown constant bias as 
\[
\xi_{\rm m} = \xi + b
\]
where $b \in \mathfrak{g}$ is the constant unknown bias. Suppose also that we measure landmarks and vectors such that an $n\times n$ matrix-valued signal $A$ of the form
\[
A = F g  
\]
or
\[
A = g^{-1}F
\] 
is available, where $F$ is an $n\times n$ invertible matrix. In this paper, we design continuous observers that globally and exponentially estimate $(g,b)$ with $\xi_{\rm m}$ and $A$, where it is assumed that the value of $F$ is available.

Relevant works are  listed in the following. In \cite{KhTrMaLa15}, the authors proposed continuous observers that estimate $(g,b)$ with $\xi_m$ and homogeneous outputs. Their observers are uniformly locally exponentially stable, but not globally exponentially stable. A similar work was done in \cite{KhTrMaLa13}, where a gradient-like innovation term was used in the observer design. The observers therein are not globally exponentially stable but only uniformly locally exponentially stable.  Gradient-like observers were also proposed in \cite{LaTrMa10}, but these observers are not globally exponentially convergent either.  

To the best of our knowledge, our observers in the present paper are the first  globally exponentially convergent {\it continuous} observers for velocity bias and state for kinematic systems on matrix Lie groups. One noticeable difference between the observers in \cite{LaTrMa10, KhTrMaLa13,KhTrMaLa15} and ours is that our observers are designed in $\mathbb R^{n\times n} \times \mathfrak{g}$ instead of $G \times \mathfrak{g}$, where $G \subset \mathbb R^{n\times n}$, such that the Euclidean structure of $\mathbb R^{n\times n}$ is fully utilized without being constrained to the group structure of $G$. This type of observers built in Euclidean space is called geometry-free and they have been widely used  for $\operatorname{SO}(3)$, e.g. \cite{BaSiOl12,MaSa17}.  

The paper is organized as follows. In Section \ref{section:main:results} we propose various forms of  globally exponentially convergent continuous observers  for velocity bias and state for kinematic systems on matrix Lie groups.  In Section \ref{section:example:SE3}, we illustrate one of the observers proposed in Section \ref{section:main:results} by applying it to the special Euclidean group $\operatorname{SE}(3)$.  The paper is concluded in Section \ref{section:conclusion}.


\section{Main Results}
\label{section:main:results}
Let $G$ be a matrix Lie group that is a subgroup of $\GLn = \{ A \in \mathbb R^{n\times n} \mid \det A \neq 0\}$, and let $\Lalge$ denote the Lie algebra of $G$. Since $G$ is a subgroup fo $\GLn$, we may assume that $\Lalge$ is a subalgebra of $(\mathbb R^{n\times n}, [\, , ])$, where $[\, ,]$ is the usual matrix commutator defined by $[A,B] = AB - BA$ for all $A,B \in \mathbb R^{n\times n}$. Let $\proj : \mathbb R^{n\times n} \rightarrow \Lalge$ denote the orthogonal projection onto $\Lalge$ with respect to the Euclidean inner product $\langle\, , \rangle$ that is defined by $\langle A, B\rangle  = \operatorname{tr}(A^TB)$ for $A,B \in \mathbb R^{n\times n}$. Let $\| \|$ denote the Euclidean or Frobenius norm which is defined by $\|A\| = \sqrt{\langle A, A\rangle}$ for all $A\in \mathbb R^{n\times n}$. For a square matrix $A$, $\lambda_{\rm minx}(A)$ and $\lambda_{\rm max}(A)$ denote the minimum eigenvalue and the maximum eigenvalue of $A$, respectively. For any matrix $A$, $\sigma_{\rm min}(A)$ and $\sigma_{\rm max}(A)$ denote the minimum singular value and the maximum singular value of $A$, respectively.
 For any $A \in \mathbb R^{n\times n}$, $\|A\|^2 = \sum_{i=1}^{n}\sigma_i^2(A)$,  where $\sigma_i(A)$'s are the singular values of $A$. 
We have $\lambda_{\rm min}(A^TA)\|B\|^2 \leq \langle AB,AB\rangle \leq \lambda_{\rm max}(A^TA)\|B\|^2$ for all $A \in \mathbb R^{n\times m}$ and $B\in \mathbb R^{m \times \ell}$, i.e.   $\sigma^2_{\rm min}(A)\|B\| \leq \| AB \| \leq \sigma^2_{\rm max}(A)\|B\|$ for all $A \in \mathbb R^{n\times m}$ and $B\in \mathbb R^{m \times \ell}$. Refer to \cite{Bl03} for more about Lie groups in the context of geometric control and mechanics. 


\subsection{Observer I}
The invariant kinematic equation on a matrix Lie group $G \subset \mathbb R^{n\times n}$ is given by
\begin{equation}\label{rigid:body}
\dot g = g \xi,
\end{equation}
where $g \in G$ and $\xi \in \Lalge$. 
Suppose that there is given an arbitrary trajectory of the system $(g(t), \xi (t)) \in G \times \Lalge$, $0 \leq t < \infty$.  We make the following three assumptions.
\begin{assumption}\label{assumption:A}
A matrix-valued signal $A (t) \in \mathbb R^{n\times n}$ is available that  can be expressed as
\begin{equation}\label{relation:AG}
A =  \meas g,
\end{equation}
 where $\meas$ is a constant  invertible matrix in $\mathbb R^{n\times n}$ and $g \in G$.
\end{assumption}

\begin{assumption}\label{assumption:Omega}
A $\Lalge$-valued signal $\xi_{\rm m}(t)$  with bias is available and related to  the true  $\xi (t) \in \Lalge$ of as follows:
\begin{align*}
 \xi_{\rm m} = \xi + b,
\end{align*}
where $b \in \Lalge$ is an unknown constant bias vector.   \end{assumption}

\begin{assumption}\label{assumption:Ombbound}
There are known constants $B_\xi>0$ and $B_b>0$ such that  $\| \xi (t) \| \leq B_\xi$ for all $t\geq 0$ and $\|b\|\leq B_b$.  There are numbers $L_g>0$ and $U_g>0$ such that
  \[
  L_g \leq \sigma_{\rm min}(g(t)) \leq \sigma_{\rm max}(g(t)) \leq U_g 
  \]
 for all $t\geq 0$, where the knowledge on the values of $L_g$ and $U_g$ is not assumed.  
\end{assumption}

We propose the following  observer:
\begin{subequations}\label{observer:Chang}
\begin{align}
\dot{\bar A} &= \bar A {\xi}_{\rm m} - A \bar b +k_P  (A - \bar A),\\
\dot{\bar b} &= - k_I \proj (A^T (A-\bar A)) \label{observer:Chang:b}
\end{align}
\end{subequations}
with  $k_P>(B_\xi + B_b)$ and $k_I >0$, where  $(\bar A, \bar b) \in \mathbb R^{n\times n} \times \Lalge$ is an estimate of $(A,b) \in G \times \Lalge$. So,  $(\meas^{-1}\bar A, \bar b)  \in \mathbb R^{n\times n} \times \Lalge$ becomes an estimate of $(g,b)  \in G \times \Lalge$ by Assumption \ref{assumption:A}. The global and exponentially convergent property of this observer is proven in the following theorem. 

\begin{theorem} \label{main:theorem:observer}
Let
\[
E_A = A-\bar A, \quad e_b = b - \bar b.
\]
Under Assumptions \ref{assumption:A} -- \ref{assumption:Ombbound}, for any $k_P> (B_\xi +B_b)$ and $k_I>0$
there exist numbers $a>0$ and $C>0$ such that
\begin{equation}\label{final:eqn}
\| E_A(t)\| + \|e_b(t)\| \leq C( \| E_A(0)\| + \|e_b(0)\| )e^{-at}
\end{equation}
for all $t\geq0$ and all $(\bar A(0), \bar b(0)) \in \mathbb R^{n\times n} \times \Lalge$. 

\begin{proof}
See Appendix. 
\end{proof}
\end{theorem}

\begin{corollary} \label{main:corollary:observer}
Suppose that Assumptions \ref{assumption:A} -- \ref{assumption:Ombbound} hold, and  let
\[
E_g = g-\meas^{-1}\bar A, \quad e_b = b - \bar b.
\]
 Then, 
there exist numbers $a>0$ and $C>0$ such that
\begin{equation}\label{final:eqn:2}
\| E_g(t)\| + \|e_b(t)\| \leq C( \| E_g(0)\| + \|e_b(0)\| )e^{-at}
\end{equation}
for all $t\geq0$ and all $(\bar A(0), \bar b(0)) \in\mathbb R^{n\times n} \times \Lalge$.
\begin{proof}
Use $\|E_g\| / \|\meas^{-1}\| \leq\|E_A\| \leq \|\meas\| \|E_g\|$ and  \eqref{final:eqn} with the constant $C$ redefined appropriately.
\end{proof}
\end{corollary}
Namely, the estimate $(\meas^{-1}\bar A(t), \bar b(t))$ converges  globally and exponentially  to the true value  $(g(t), b)$ as $t$ tends to $\infty$.

\begin{remark}
We can also build an observer that allows $\bar b$ to be in $\mathbb R^{n\times n}$ instead of $\mathfrak{g}$. The modified observer is given by
\begin{subequations}\label{observer:Chang:mod}
\begin{align}
\dot{\bar A} &= \bar A {\xi}_{\rm m} - A \bar b +k_P  (A - \bar A),\\
\dot{\bar b} &= - k_I A^T (A-\bar A),  \label{observer:Chang:mod:b}
\end{align}
\end{subequations}
where $(\bar A, \bar b) \in \mathbb R^{n\times n} \times \mathbb R^{n\times n}$. 
Notice that the projection operator $\proj$ in \eqref{observer:Chang:b} is removed from \eqref{observer:Chang:b} to obtain \eqref{observer:Chang:mod:b}. Theorem \ref{main:theorem:observer} and Corollary \ref{main:corollary:observer} also hold for this observer, whose proof is almost identical to the proofs of Theorem \ref{main:theorem:observer} and Corollary \ref{main:corollary:observer}, so it is left to the reader.
\end{remark}

\begin{remark}
   Assumption \ref{assumption:A} can be relaxed by allowing the matrix $F$ to be time-varying. More specifically, we make the following assumption: there  are numbers $\ell_{\rm min} >0$  and $\ell_{\rm max} >0$ such that 
 \begin{equation}\label{relax:G}
 \ell_{\rm min} \leq  \sigma_{\rm min} (\meas(t)) \leq \sigma_{\rm max} (\meas (t)) \leq \ell_{\rm max}
 \end{equation}
  for all $t \geq 0$. 
In this case, we propose the following  observer:
\begin{align*}
\dot{\bar A} &= \bar A {\xi}_{\rm m} - A {\bar b} +k_P  (A - \bar A) + \dot \meas \meas^{-1}A,\\
\dot{\bar b} &= -k_I \proj (A^T(A- \bar A))
\end{align*}
with  $k_P>0$ and $k_I >0$, where  $(\bar A, \bar b) \in \mathbb R^{n\times n} \times \Lalge$ is an estimate of $(A,b)$.  It is not difficult to show that Theorem  \ref{main:theorem:observer}  and Corollary \ref{main:corollary:observer} also hold for this observer with the relaxed assumption on $\meas(t)$ as above.  Here, the knowledge on the values of $\ell_{\rm min} $  and $\ell_{\rm max}$ is not required here. 
\end{remark}

\begin{remark}
Since the estimate $F^{-1}\bar A$ may not lie in $G$ in general, one may need to project it to $G$ as an output of the observer although $F^{-1}\bar A(t)$ converges to $g(t) \in G$ as $t$ tends to infity. For example, if $G = \SO$, then the usual polar decomposition can be used to define a projection from $\mathbb R^{3\times 3}$ to $\SO$. Projection for  $\operatorname{SE}(3)$ will be discussed in Section \ref{section:example:SE3}. However, if one designs controllers in $\mathbb R^{n\times n}$ for an extension of \eqref{rigid:body} into $\mathbb R^{n\times n}$ as proposed in \cite{Ch18b}, then the direct use of $F^{-1}\bar A$ in feedback would be fine. 
\end{remark}


\subsection{Observer II}

Recall the kinematic equation in \eqref{rigid:body}.
 We now consider a case where the measurement matrix $A$ is related to the true signal $g(t)$ as $A = g^{-1}(t)\meas$ instead of $A = \meas g(t)$. Consequently,  in place of Assumption \ref{assumption:A:beta}, let us make the following assumption:
\begin{assumption}\label{assumption:A:beta}
A matrix-valued signal $A (t) \in \mathbb R^{n\times n}$ is available that  can be expressed as
\begin{equation}\label{relation:AG:beta}
A =  g^{-1}\meas,
\end{equation}
 where $\meas$ is a constant  invertible matrix in $\mathbb R^{n\times n}$ and $g \in G$.
\end{assumption}
By \eqref{rigid:body}, $A$ defined in \eqref{relation:AG:beta} satisfies
\begin{equation}\label{Adot:beta}
\dot A = -\xi A.
\end{equation}
Under Assumptions \ref{assumption:A:beta}, \ref{assumption:Omega} and \ref{assumption:Ombbound}, 
we propose the following  observer:
\begin{subequations}\label{observer:Chang:beta}
\begin{align}
\dot{\bar A} &= - {\xi}_{\rm m} \bar A + \bar b A +k_P  (A - \bar A), \label{observer:Chang:beta:A}\\
\dot{\bar b} &=  k_I \proj ( (A-\bar A)A^T) \label{observer:Chang:beta:b}
\end{align}
\end{subequations}
with  $k_P>(B_\xi + B_b)$ and $k_I >0$, where  $(\bar A, \bar b) \in \mathbb R^{n\times n} \times \Lalge$ is an estimate of $(A,b) \in G \times \Lalge$.
\begin{theorem} \label{main:theorem:observer:beta}
For the observer \eqref{observer:Chang:beta},
let
\[
E_A = A-\bar A, \quad e_b = b - \bar b.
\]
Under  Assumptions \ref{assumption:A:beta}, \ref{assumption:Omega} and \ref{assumption:Ombbound}, for any $k_P> (B_\xi + B_b)$ and $k_I>0$
there exist numbers $a>0$ and $C>0$ such that
\[
\| E_A(t)\| + \|e_b(t)\| \leq C( \| E_A(0)\| + \|e_b(0)\| )e^{-at}
\]
for all $t\geq0$ and all $(\bar A(0), \bar b(0)) \in \mathbb R^{n\times n} \times \Lalge$. 
\begin{proof}
See Appendix.
\end{proof}
\end{theorem}

\begin{corollary} \label{main:corollary:observer:beta}
Consider the observer \eqref{observer:Chang:beta}.
Suppose that Assumptions \ref{assumption:A:beta}, \ref{assumption:Omega} and \ref{assumption:Ombbound} hold, and  let
\[
E_g = g-\meas \bar A^{-1}, \quad e_b = b - \bar b.
\]
 Then, 
there exist numbers $a>0$ and $C>0$ such that
\[
\| E_g(t)\| + \|e_b(t)\| \leq C( \| E_g(0)\| + \|e_b(0)\| )e^{-at}
\]
for all $t\geq0$ and all $(\bar A(0), \bar b(0)) \in\mathbb R^{n\times n} \times \Lalge$.
\end{corollary}
In other words, the estimate $(\meas \bar A(t)^{-1}, \bar b(t))$ converges  globally and exponentially  to the true value  $(g(t), b)$ as $t$ tends to infinity. 
\begin{remark}
If $\meas (t)$ is time-varying such that \eqref{relax:G} is satisfied for all $t\geq0$, then \eqref{observer:Chang:beta:A} has only to be modified to
\[
\dot{\bar A} = - {\xi}_{\rm m} \bar A + \bar b A +k_P  (A - \bar A) + A \meas^{-1}\dot \meas
\]
while \eqref{observer:Chang:beta:b} remains intact. 
\end{remark}

 We now derive from \eqref{observer:Chang} various  observers of {\it concrete} form that estimate $(R,b)$ from vector measurements. 
  Assume that there is a set $\mathcal S = \{s_i, 1\leq i \leq m\}$ of $m$ known fixed inertial vectors, where each $s_i$ in $\mathcal S$ is a vector in $\mathbb R^n$, such that the rank of $\mathcal S$ is $n$.  Assume also that measurements of the vectors are made in the body-fixed frame and the set of the measured vectors is denoted by $\mathcal C = \{ c_i , 1\leq i \leq m\}$ and  related to $\mathcal S$ as follows:
 \[
 c_i = g^{-1} s_i, \quad i = 1, \ldots, m,
 \]
 where $\meas \in G$.  Let
 \begin{equation}\label{def:SC}
 S =  \begin{bmatrix}
 s_1& \cdots & s_m
 \end{bmatrix}, \quad C =  \begin{bmatrix}
 c_1& \cdots & c_m
 \end{bmatrix}
 \end{equation}
  be $n \times m$ matrices made of the {\it column} vectors from $\mathcal S$ and ${\mathcal C}$, respectively.

 \begin{corollary}\label{corollary:obs:linear:beta}
Let  $S$ and $C$ be given  in \eqref{def:SC}. If there is a matrix $W \in \mathbb R^{m\times n}$ such that $\meas:= SW$ has rank $n$, then \eqref{relation:AG:beta} is satisfied by $A = CW$ and  the observer \eqref{observer:Chang:beta} is applicable.
 \begin{proof}
Trivial.
 \end{proof}
 \end{corollary}

\begin{remark}\label{remark:obs:linear:beta}
 Corollary  \ref{corollary:obs:linear:beta} can be applied in several ways. For example, 
the substitution of $WS^T$ into $W$ in Corollary \ref{corollary:obs:linear:beta} would yield
\[
 \meas = SWS^T, \quad   A = CWS^T,
\]
 where it is assumed that $W$ is an $m\times m$ matrix  such that $\meas$ has rank $n$. Likewise,  $W$ in Corollary  \ref{corollary:obs:linear:beta} can be chosen such that $F$ depends more nonlinearly on $S$.
 \end{remark}


\subsection{Variants}
We here propose an observer that is a variant of the observer \eqref{observer:Chang} with $A^T$ replaced by $A^{-1}$ in \eqref{observer:Chang:b}.  Recall the kinematic equation \eqref{rigid:body},  and under Assumptions \ref{assumption:A} -- \ref{assumption:Ombbound},
we propose the following  new observer:
\begin{subequations}\label{observer:Chang:three}
\begin{align}
\dot{\bar A} &= \bar A {\xi}_{\rm m} - A \bar b +k_P  (A - \bar A),\\
\dot{\bar b} &= - k_I \proj (A^{-1} (A-\bar A)) 
\end{align}
\end{subequations}
with  $k_P> (2B_\xi + B_b)$ and $k_I >0$, where  $(\bar A, \bar b) \in \mathbb R^{n\times n} \times \Lalge$ is an estimate of $(A,b) \in G \times \Lalge$. So,  $(\meas^{-1}\bar A, \bar b)  \in \mathbb R^{n\times n} \times \Lalge$ becomes an estimate of $(g,b)  \in G \times \Lalge$ by Assumption \ref{assumption:A}.  

\begin{theorem} \label{main:theorem:observer:three}
For the observer \eqref{observer:Chang:three},
let
\[
E_A = A-\bar A, \quad e_b = b - \bar b.
\]
Under Assumptions \ref{assumption:A} -- \ref{assumption:Ombbound}, for any $k_P> (2 B_\xi +B_b)$ and $k_I>0$
there exist numbers $a>0$ and $C>0$ such that
\[
\| E_A(t)\| + \|e_b(t)\| \leq C( \| E_A(0)\| + \|e_b(0)\| )e^{-at}
\]
for all $t\geq0$ and all $(\bar A(0), \bar b(0)) \in \mathbb R^{n\times n} \times \Lalge$. 

\begin{proof}
See Appendix. 
\end{proof}
\end{theorem}

We also  propose a variant of the observer \eqref{observer:Chang:beta} with $A^T$ replaced by $A^{-1}$ in \eqref{observer:Chang:beta:b}.  Under Assumptions \ref{assumption:Omega},  \ref{assumption:Ombbound} and \ref{assumption:A:beta}, we propose the following new observer:
\begin{subequations}\label{observer:Chang:four}
\begin{align}
\dot{\bar A} &= - {\xi}_{\rm m} \bar A  +\bar b A  +k_P  (A - \bar A),\\
\dot{\bar b} &= k_I \proj ( (A-\bar A)A^{-1}) 
\end{align}
\end{subequations}
with  $k_P> (2B_\xi + B_b)$ and $k_I >0$, where  $(\bar A, \bar b) \in \mathbb R^{n\times n} \times \Lalge$ is an estimate of $(A,b) \in G \times \Lalge$. So,  $(\meas \bar A^{-1}, \bar b)  \in \mathbb R^{n\times n} \times \Lalge$ becomes an estimate of $(g,b)  \in G \times \Lalge$ by Assumption \ref{assumption:A}.   
\begin{theorem} 
For the observer \eqref{observer:Chang:four},
let
\[
E_A = A-\bar A, \quad e_b = b - \bar b.
\]
Under  Assumptions \ref{assumption:A:beta}, \ref{assumption:Omega} and \ref{assumption:Ombbound}, for any $k_P> (B_\xi + B_b)$ and $k_I>0$
there exist numbers $a>0$ and $C>0$ such that
\[
\| E_A(t)\| + \|e_b(t)\| \leq C( \| E_A(0)\| + \|e_b(0)\| )e^{-at}
\]
for all $t\geq0$ and all $(\bar A(0), \bar b(0)) \in \mathbb R^{n\times n} \times \Lalge$. 
\begin{proof}
Omitted since it is similar to the proof of Theorem \ref{main:theorem:observer:three}.
\end{proof}
\end{theorem}

\begin{remark}
1. Corollaries \ref{main:corollary:observer} and  \ref{main:corollary:observer:beta}  also hold for the observers \eqref{observer:Chang:three} and \eqref{observer:Chang:four}, respectively.

2. Corollary \ref{corollary:obs:linear:beta} and Remark \ref{remark:obs:linear:beta} also hold true for the observer \eqref{observer:Chang:four}.
\end{remark}


\section{Example: Application to $\operatorname{SE}(3)$}
\label{section:example:SE3}

\begin{figure}
\begin{center}
\includegraphics[scale = 0.23]{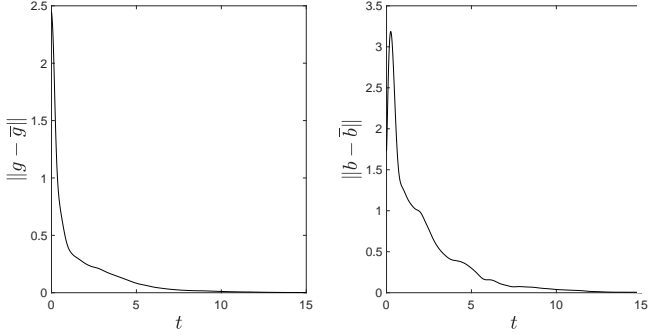}
\end{center}
\caption{\label{figure:simulation_SE3} The pose estimation error $\| g(t) - \bar g(t)\|$  and the velocity bias estimation error $\|b - \bar b(t)\|$ by the observer \eqref{observer:Chang:beta}. }
\end{figure}

 We now illustrate the theory presented in Section \ref{section:main:results} with the special Euclidean  group on $\mathbb R^3$. The group can be expressed in homogeneous coordinates as 
\[
\SE = \left \{ \begin{bmatrix}
R & x \\ 0 & 1
\end{bmatrix} \mid R \in \SO, x \in \mathbb R^3 \right \},
\]
where $\SO = \{ R\in \mathbb R^{3\times 3}\mid R^T R = I, \det R = 1\}$ is the special orthogonal group whose Lie algebra is $\mathfrak{so}(3) = \{ A \in \mathbb R^{3 \times 3} \mid A ^T = -A\}$. 
It is easy to see that $\SE$ is a subgroup of $\operatorname{GL}(4) = \{ A \in \mathbb R^{4\times 4} \mid \det A \neq 0\}$. The Lie algebra of $\SE$  is then given by
\[
\mathfrak{se}(3) =\left \{\begin{bmatrix}
\hat\Omega & v \\ 0 & 0
\end{bmatrix} \textup{ for some } \Omega \in \mathbb R^{3}, x \in \mathbb R^3 \right \},
\]
where the hat map $\wedge : \mathbb R^3 \rightarrow \mathfrak{so}(3)$ is defined such that $\hat x y = x \times y$ for all $x,y \in \mathbb R^3$. 
In homogeneous coordinates, landmarks to be measured with sensors are expressed  in the form  
\begin{equation}
\begin{bmatrix}
x \\ 1
\end{bmatrix}, \quad x \in \mathbb R^3,
\end{equation}
and such vectors {\it at infinity}  as the gravity or the Earth's magnetic field are expressed  in the form
\begin{equation}
\begin{bmatrix}
x \\ 0
\end{bmatrix}, \quad x \in \mathbb R^3.
\end{equation} 
The orthogonal projection $\pi_{\mathfrak{se}(3)}: \mathbb R^{4\times 4} \rightarrow \mathfrak{ se}(3)$ is given as follows: For $A \in \mathbb R^{4\times 4}$ given by
\[
A = \begin{bmatrix} 
B & x \\ y^T & z
\end{bmatrix}, \quad B \in \mathbb R^{3\times 3}, x \in \mathbb R^{3\times 1},  y \in \mathbb R^{3\times 1},z\in \mathbb R,
\]
we have
\[
\pi_{\mathfrak{se}(3)}(A) = \begin{bmatrix}
\frac{1}{2}(B - B^T) & x \\
0_{1\times 3} & 0
\end{bmatrix} \in \mathfrak{se}(3).
\]
Suppose that we measure  in the body frame the following inertial vectors given by
\begin{align*}
&s_1 = (e_1,1), s_2 = (e_2,1), s_3 = ( e_3,1), s_4 = (e_1 + e_3,1), \\
&s_5 = (-e_3,0),
\end{align*}
where $\{e_1,e_2,e_3\}$ is the standard basis of $\mathbb R^3$ and $s_5$ represents the gravity direction.  Suppose the measured signal matrix $A(t)$ is given by
\begin{align*}
A(t) = g(t)^{-1}\meas =  g(t)^{-1}SWS^T = C(t) WS^T 
\end{align*}
with $\meas = SWS^T$ and $C(t) = g(t)^{-1}S$, where
\[
S = \begin{bmatrix} 
 s_1 & s_2 & s_3 & s_4 & s_5
\end{bmatrix}.
\]
Here, each column in the $C(t)$ matrix is what is measured in the body-fixed frame. 
For convenience, we set $W = I_{4\times 4}$ although any $4\times 4$  matrix such that $SWS^T$ is invertible would work for $W$. 
Suppose that a set of  true trajectories $(R(t), x(t)) \in \SE$ and $(\Omega (t), V(t)) \in \mathbb R^3 \times \mathbb R^3$ are given as follows: 
\begin{align}
R(t)&=\exp(t\hat e_1)\exp(t\hat e_3)\exp(t\hat e_1),  \label{true:R} \\
x(t) &= (\cos t, \sin t, \cos t),\\
 \Omega(t)&=(1+\cos t, \sin t - \sin t\cos t, \cos t+\sin^2 t), \label{true:Omega}\\
 V(t) &= R^T(t)\dot x(t),
\end{align}
where $\Omega(t)$ satisfies $\hat \Omega (t) = R^T(t)\dot R(t)$. 
  Assume  that the  unknown constant gyro bias $b_\Omega$ and the unknown constant velocity bias $b_v$ are respectively  given by 
 \begin{equation}\label{true:bias}
 b_\Omega=(1, 0.5,-1), \quad b _v = (0.5,-0.5,0.5).
 \end{equation}  
We use the observer of the form \eqref{observer:Chang:beta}. The  gains are chosen as $k_P = 4$ and $k_I = 0.75$,  and the initial state of the observer is given by 
\[
\bar A (0) = \bar g_0^{-1} 
\meas
\]
where
\[
\bar g_0 = \begin{bmatrix}
\exp(\displaystyle{\frac{\pi}{2}\hat e_1}) & 0_{3\times 1}\\
0_{1\times 3} & 1
\end{bmatrix}
\]
 and 
 \[
 \bar b_\Omega (0) = (0,0,0), \quad \bar b_v(0) = (0,0,0).
 \]
The simulation results are plotted in Fig. \ref{figure:simulation_SE3}, where the pose estimation error $\| g(t) - \bar g(t)\|$
with  $\bar g(t) := F \bar A(t)^{-1}\in \mathbb R^{4\times 4}$,
and the bias estimation error $\|b - \bar b(t)\|$ are plotted. It can be seen that both estimation errors converge well to zero as theoretically predicted.

\begin{figure}
\begin{center}
\includegraphics[scale = 0.23]{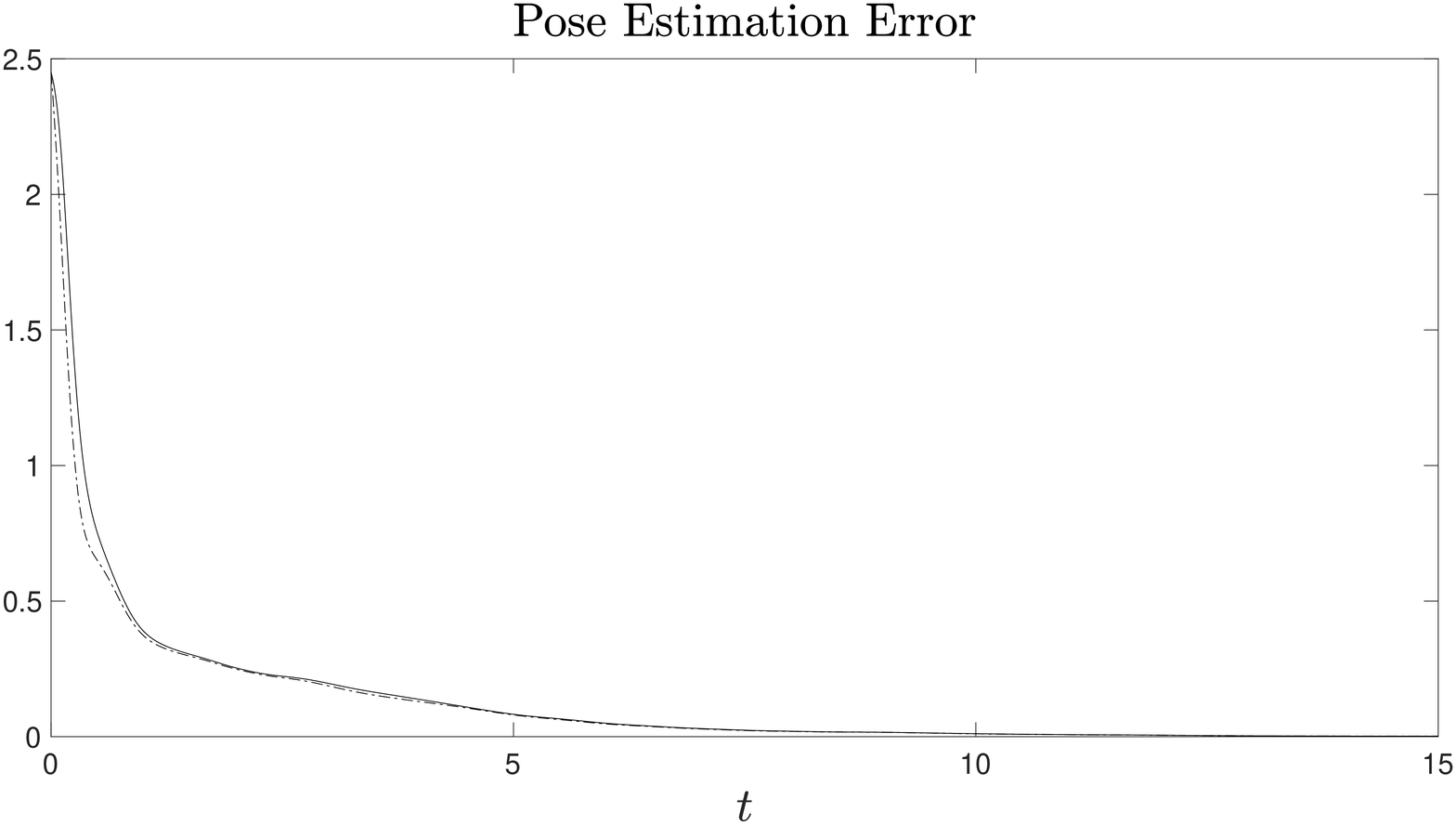}
\end{center}
\caption{\label{figure:SE3_errors} Two pose estimation errors by the observer \eqref{observer:Chang:beta}: $\| g(t) - \bar g(t)\|$ by the  observer  (solid) and  $\| g(t) - \bar g_{\operatorname{SE}(3)}(t)\|$ by the $\operatorname{SE}(3)$ factor $\bar g_{\operatorname{SE}(3)}(t)$ of  $\bar g(t)$  obtained through the projection \eqref{SE3:proj}  (dash-dot).}
\end{figure}

To examine if the image trajectory of $\bar g(t)$ under a projection onto $\operatorname{SE}(3)$  also converges to $g(t)$, let us 
define a projection $\operatorname{proj} : \mathbb R^{4\times 4} \rightarrow \operatorname{SE}(3)$ 
as follows: for  any 
\[
\bar g = \begin{bmatrix}
\bar g_1 & \bar g_2 \\
\bar g_3 & \bar g_4
\end{bmatrix} \in \mathbb R^{4\times 4}
\]
with $\bar g_1 \in \mathbb R^{3\times 3}$, $\bar g_2 \in \mathbb R^{3\times 1}$,  $\bar g_3 \in \mathbb R^{1\times 3}$, and $\bar g_4\in \mathbb R$,
\begin{equation}\label{SE3:proj}
\operatorname{proj} (\bar g) := \begin{bmatrix}
\bar g_{1,\SO} & \bar g_2 \\
0_{1\times 3} & 1
\end{bmatrix} \in \operatorname{SE}(3),
\end{equation}
where $\bar g_{1,\SO}$ denotes the $\SO$ factor  in  polar decomposition of  $\bar g_1$. For convenience, let
\[
\bar g_{\operatorname{SE}(3)} (t):= \operatorname{proj} (\bar g (t)).
\]
The pose estimation error $\| g (t) - \bar g_{\operatorname{SE}(3)} (t)\|$ by $\bar g_{\operatorname{SE}(3)} (t)$ is plotted in Fig. \ref{figure:SE3_errors} along with the pose estimation error $\| g(t) - \bar g(t)\|$ by $\bar g(t)$ that was obtained in the simulation. It can be seen in the figure that $\bar g(t)$ stays very close to its $\operatorname{SE}(3)$ factor $ \bar g_{\operatorname{SE}(3)} (t)$, and $ \bar g_{\operatorname{SE}(3)} (t)$ also converges to the true pose $g(t)$ as time tends to infinity.

\begin{figure}
\begin{center}
\includegraphics[scale = 0.23]{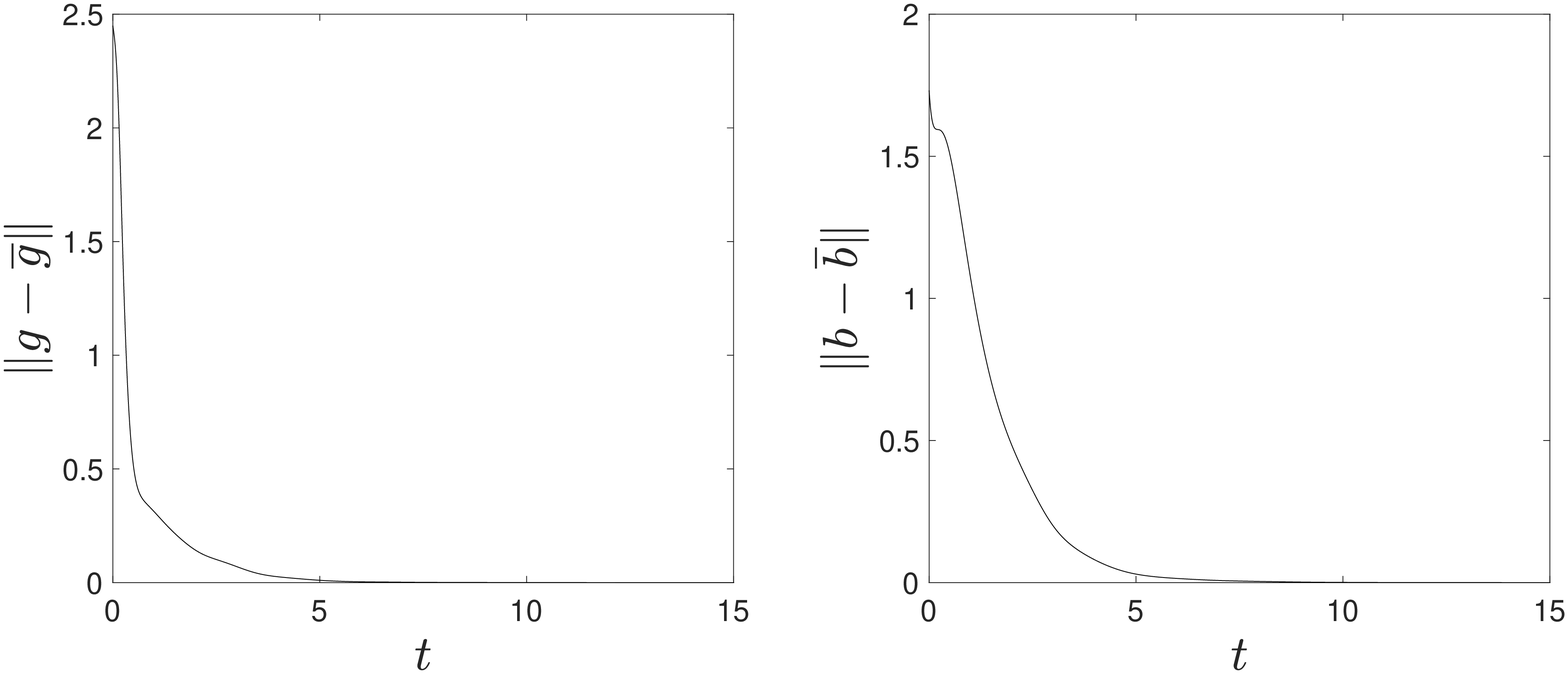}
\end{center}
\caption{\label{figure:simulation_SE3_2} The pose estimation error $\| g(t) - \bar g(t)\|$  and the velocity bias estimation error $\|b - \bar b(t)\|$ by the observer \eqref{observer:Chang:four}. }
\end{figure}

\begin{figure}
\begin{center}
\includegraphics[scale = 0.23]{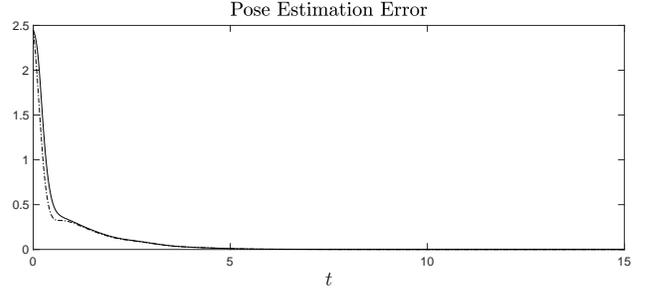}
\end{center}
\caption{\label{figure:SE3_errors_2} Two pose estimation errors by the observer \eqref{observer:Chang:four}: $\| g(t) - \bar g(t)\|$ by the  observer  (solid) and  $\| g(t) - \bar g_{\operatorname{SE}(3)}(t)\|$ by the $\operatorname{SE}(3)$ factor $\bar g_{\operatorname{SE}(3)}(t)$ of  $\bar g(t)$  obtained through the projection \eqref{SE3:proj}  (dash-dot).}
\end{figure}

For the purpose of comparison, we now apply the observer \eqref{observer:Chang:four} with the same setting except the observer gains which are now chosen as $k_P = 4$ and $k_I = 4$. The estimation results are plotted in Fig. \ref{figure:simulation_SE3_2}. It can be seen  that the bias estimation error by the observer \eqref{observer:Chang:four} in Fig.  \ref{figure:simulation_SE3_2} converges fast without overshoot in comparison with the estimation error  by \eqref{observer:Chang:beta} that is plotted in Fig.  \ref{figure:simulation_SE3}. The $\operatorname{SE}(3)$ part of $\bar g(t)$ computed by the projection \eqref{SE3:proj} also converges well to the true signal $g(t)$ as shown in Fig. \ref{figure:SE3_errors_2}.

\begin{remark}
There have been papers on estimation of pose and velocity measurement bias for $\operatorname{SE}(3)$, e.g. \cite{WaTa18,VaCuSiOl10,HuHaMaTr15} and references therein.  A globally exponentially convergent hybrid  (not continuous) observer is proposed in \cite{WaTa18}, and a non-global exponentially convergent observer is proposed in \cite{VaCuSiOl10}.  A gradient-like observer design on $\operatorname{SE}(3)$ with system outputs on the real projective space was proposed in \cite{HuHaMaTr15}. Refer to \cite{GoLe17} for a global formulation of  extended Kalman filter on $\operatorname{SE}(3)$ for geometric control of a drone. 
\end{remark}

\section{Conclusion}
\label{section:conclusion}
We have successfully designed   globally exponentially convergent continuous observers for  kinematic  invariant systems on finite-dimensional matrix Lie groups that   
estimate state and constant velocity bias from measurements of landmarks, vectors and biased velocity. We have applied the result to the special Euclidean group $\operatorname{SE}(3)$ and carried out a simulation study to illustrate an excellent performance of the observer for $\operatorname{SE}(3)$.
We  plan to apply the result to drone control \cite{ChEu17} and to combine it with deep neural networks \cite{CaCh18}.

\appendix

\subsection*{ Proof of Theorem \ref{main:theorem:observer}}
\begin{proof}
From \eqref{rigid:body} and Assumption \ref{assumption:A}, $A(t)$ satisfies
\begin{equation}\label{Adot:eq}
\dot A =  A \xi. 
\end{equation}
By Assumption \ref{assumption:Omega}, the observer \eqref{observer:Chang} can be written as
\begin{subequations}\label{modified:observer:Chang}
\begin{align}
\dot{\bar A} &= \bar A(\xi +  b) - A \bar b +k_P E_A,\\
\dot{\bar b} &= -k_I \proj (A^T E_A).
\end{align}
\end{subequations}
  By Assumption \ref{assumption:Ombbound}, there is a number $\epsilon$ such that
\[
0 < \epsilon < \min \left \{ H, \frac{1}{\|\meas \| U_g\sqrt{k_I}} \right \},
\]
where
\[
H=  \frac{4(k_P-B_\xi- B_b) L _g^2\lambda_{\min} (\meas^T\meas) }{ (4k_I L_g^2 \lambda_{\rm min}(\meas^T\meas) +  (k_P + B_b + 2B_\xi)^2 )U_g^2 \|\meas\|^2 }.
\]
The following three quadratic functions of $(\|E_A\|, \|e_b\|)$ are then all positive definite:
\begin{align*}
V_1(\|E_A\|, \|e_b\|) &= \frac{1}{2} \|E_A\|^2 + \frac{1}{2k_I}\|e_b\|^2 \! - \! \epsilon U_g\|\meas\| \|E_A\| \|e_b\|,\\
V_2(\|E_A\|, \|e_b\|) &= \frac{1}{2} \|E_A\|^2 + \frac{1}{2k_I}\|e_b\|^2 \! + \! \epsilon U_g\|\meas\| \|E_A\| \|e_b\|,\\
V_3(\|E_A\|, \|e_b\|) &= (k_P -(B_\xi + B_b)- \epsilon k_IU_g^2 \|\meas\|^2) \|E_A\|^2 \\
&\quad+ \epsilon  \lambda_{\rm min}(\meas^T \meas)L_{g}^2\|e_b\|^2 \\
&\quad -\epsilon (k_P + B_b+2B_\xi)U_g\|\meas\|\|E_A\| \|e_b\|.
\end{align*}
Hence, there are numbers $\alpha>0$ and $\beta >0$ such that 
\begin{equation}\label{V213}
V_2 \leq \alpha V_1, \quad \beta V_2 \leq  V_3.
\end{equation}
Let 
\[
V(E_A, e_b) = \frac{1}{2} \|E_A\|^2 + \frac{1}{2k_I} \|e_b\|^2 + \epsilon \langle E_A, A e_b \rangle,
\]
which satisfies
\begin{equation}\label{V1:V:V2}
V_1(\|E_A\|, \|e_b\|) \leq V(E_A, e_b) \leq V_2(\|E_A\|, \|e_b\|)
\end{equation}
for all $(E_A, e_b) \in \mathbb R^{n\times n} \times \Lalge$ by the Cauchy-Schwarz inequality and $\|A\| = \|\meas g\| \leq \|\meas\| U_g$. From  \eqref{Adot:eq}, \eqref{modified:observer:Chang},  and the assumption of the bias $b$ being constant, it follows that the estimation error $(E_A, e_b)$ obeys
\begin{align*}
\dot E_A &= E_A(\xi +  b) - A e_b - k_PE_A,\\
\dot e_b &= k_I \proj(A^TE_A).
\end{align*}
Along any trajectory of the composite system consisting of the rigid body \eqref{rigid:body} and the observer \eqref{observer:Chang},
\begin{align*}
\frac{dV}{dt} &= \langle E_A, E_A(\xi +  b) - A e_b - k_PE_A\rangle \\
&\quad+ \langle  e_b, \proj (A^TE_A) \rangle \\
&\quad + \epsilon \langle E_A (\xi + b) - A e_b - k_PE_A, A e_b\rangle  \\
&\quad+ \epsilon \langle E_A, A\xi e_b\rangle +\epsilon k_I \langle E_A, A\proj (A^TE_A)\rangle \\
&\leq - (k_P -(B_\xi + B_b)- \epsilon k_IU_g^2 \|\meas\|^2) \|E_A\|^2 \\
&\quad- \epsilon  \lambda_{\rm min}(\meas^T \meas)L_g^2\|e_b\|^2 \\
&\quad +\epsilon (k_P + B_b+2B_\xi)U_g\|\meas\|\|E_A\| \|e_b\| \\
&= -V_3 \leq -\beta V_2 \leq - \beta V,
\end{align*}
where the following have been used:
\begin{align*}
&\langle E_A, E_A(\xi +b) \rangle \leq \|E_A\|^2 (B_\xi + B_b),\\
&  \langle E_A, A e_b\rangle =  \langle A^TE_A,  e_b\rangle =  \langle \proj (A^TE_A),  e_b\rangle,\\
&\langle Ae_b, A e_b\rangle \geq \lambda_{\rm min}(\meas^T\meas)\|g e_b\|^2 \geq  \lambda_{\rm min}(\meas^T\meas)L_g^2\|e_b\|^2,\\
&\langle E_A, A\proj (A^TE_A)\rangle = \|\proj (A^TE_A)\|^2\\
&\qquad\qquad\qquad\qquad \quad \leq \|A^TE_A\|^2 \leq U_g^2\|\meas \|^2 \|E_A\|^2.
\end{align*}
Hence, $V(t) \leq V(0) e^{-\beta t}$ for all $t\geq 0$ and all $(\bar A(0), \bar b(0)) \in \mathbb R^{n\times n} \times \Lalge$. It follows from \eqref{V213} and \eqref{V1:V:V2} that
\[
V_1(t) \leq V(t) \leq V(0) e^{-\beta t}  \leq V_2(0) e^{-\beta t}  \leq \alpha V_1(0) e^{-\beta t} 
\]
for all $t\geq 0$ and all $(\bar A(0), \bar b(0)) \in \mathbb R^{n\times n} \times \Lalge$. Since $0<\epsilon < 1/ (\|\meas \|U_g\sqrt{k_I})$,  the map defined by
\[
(x_1, x_2) \mapsto \sqrt{\frac{1}{2}x_1^2 + \frac{1}{2 k_I}x_2^2 -  \epsilon U_g\|\meas\| x_1x_2}
\]
is a norm on $\mathbb R^2$, where $(x_1,x_2) \in \mathbb R^2$, which is equivalent to the 1-norm on $\mathbb R^2$ since all norms are equivalent on a finite-dimensional vector space. Hence,  $V_1(t) \leq \alpha V_1(0)e^{-\beta t}$ implies that there exists $C>0$ such that \eqref{final:eqn} holds  for all $t\geq 0$ and all $(\bar A(0), \bar b(0)) \in \mathbb R^{n\times n} \times \Lalge$, where $a = \beta/2$. 
\end{proof}

\subsection*{ Proof of Theorem \ref{main:theorem:observer:beta}}
\begin{proof}
  By Assumption \ref{assumption:Ombbound},  there is a number $\epsilon$ such that
\[
0 < \epsilon < \min \left \{ H, \frac{L_g}{\|\meas \| \sqrt{k_I}} \right \},
\]
where
\[
H=  \frac{4(k_P -B_\xi -B_b)L_g^2\lambda_{\min} (\meas^T\meas) }{  (4k_I  \lambda_{\rm min}(\meas^T\meas) +  (k_P + B_b + 2B_\xi)^2 U_g^2)\|\meas\|^2 }.
\]
The following three quadratic functions of $(\|E_A\|, \|e_b\|)$ are then all positive definite:
\begin{align*}
V_1(\|E_A\|, \|e_b\|) &= \frac{1}{2} \|E_A\|^2 + \frac{1}{2k_I}\|e_b\|^2 \! - \! \frac{\epsilon}{L_g} \|\meas\| \|E_A\| \|e_b\|,\\
V_2(\|E_A\|, \|e_b\|) &= \frac{1}{2} \|E_A\|^2 + \frac{1}{2k_I}\|e_b\|^2 \! + \!  \frac{\epsilon}{L_g}\|\meas\| \|E_A\| \|e_b\|,\\
V_3(\|E_A\|, \|e_b\|) &= \left (k_P -(B_\xi + B_b)- \frac{\epsilon k_I \|\meas \|^2}{L_g^2} \right ) \|E_A\|^2 \\
&\quad+ \frac{\epsilon  \lambda_{\rm min}(\meas^T \meas)}{U_{g}^2}\|e_b\|^2 \\
&\quad -\frac{\epsilon (k_P + B_b+2B_\xi)}{L_g}\|\meas\|\|E_A\| \|e_b\|.
\end{align*}
Hence, there are numbers $\alpha>0$ and $\beta >0$ such that  \eqref{V213} holds.
Let
\[
V(E_A, e_b) = \frac{1}{2} \|E_A\|^2 + \frac{1}{2k_I} \|e_b\|^2 - \epsilon \langle E_A,  e_b A \rangle,
\]
which satisfies \eqref{V1:V:V2}. Since $b$ is constant by assumption, it follows from \eqref{Adot:beta} and \eqref{observer:Chang:beta} that
\begin{align*}
\dot E_A &= -\xi_{\rm m}E_A + e_b A - k_PE_A,\\
\dot e_b &= -k_I \proj(E_A A^T).
\end{align*}
Along any trajectory of the composite system consisting of the rigid body \eqref{rigid:body} and the observer \eqref{observer:Chang:beta},
\begin{align*}
\frac{dV}{dt} &= \langle E_A, -\xi_{\rm m}E_A+e_b A  - k_PE_A\rangle - \langle  e_b, \proj (E_A A^T) \rangle \\
&\quad - \epsilon \langle -\xi_{\rm m}E_A +e_b A  - k_PE_A,  e_b A\rangle  \\
&\quad + \epsilon \langle E_A,  e_b \xi A\rangle +\epsilon k_I \langle E_A, \proj (E_A A^T)A\rangle \\
&\leq - \left (k_P -(B_\xi + B_b)- \frac{\epsilon k_I \|\meas \|^2}{L_g^2} \right ) \|E_A\|^2 \\
&\quad- \frac{ \epsilon  \lambda_{\rm min}(\meas^T \meas)}{U_g^2}\|e_b\|^2 \\
&\quad +\frac{\epsilon (k_P + B_b+2B_\xi)\|\meas\|}{L_g}\|E_A\| \|e_b\| \\
&= -V_3 \leq -\beta V_2 \leq - \beta V.
\end{align*}
The rest of the proof is identical to the corresponding part in the proof of Theorem \ref{main:theorem:observer}, so it is omitted.
\end{proof}

\subsection*{Proof of Theorem \ref{main:theorem:observer:three}}
\begin{proof}
The measured matrix $A = \meas g$ obeys \eqref{Adot:eq}
Let
\[
\mathcal{E}_A = I - A^{-1}\bar A, \quad e_b = b - \bar b.
\]
From \eqref{Adot:eq} and \eqref{observer:Chang:three},
\begin{align*}
\dot{\mathcal E}_A &=  {\mathcal E}_A\xi_{\rm m}  -\xi  {\mathcal E}_A - e_b - k_P {\mathcal E}_A\\
\dot e_b &= k_I \proj ( {\mathcal E}_A).
\end{align*}
There is an $\epsilon >0$ such that
\[
0 <\epsilon < \min \left  \{\frac{1}{\sqrt{k_I}}, \frac{4 (k_P - 2B_\xi - B_b)}{4k_I + (k_P + 2B_\xi + B_b)^2} \right  \}.
\]
The following three quadratic functions of $(\|E_A\|, \|e_b\|)$ are then all positive definite:
\begin{align*}
V_1(\|\mathcal{E}_A\|, \|e_b\|) &= \frac{1}{2} \|\mathcal{E}_A\|^2 + \frac{1}{2k_I}\|e_b\|^2 \! - \! \epsilon  \|\mathcal{E}_A\| \|e_b\|,\\
V_2(\|\mathcal{E}_A\|, \|e_b\|) &= \frac{1}{2} \|\mathcal{E}_A\|^2 + \frac{1}{2k_I}\|e_b\|^2 \! + \! \epsilon \|\mathcal{E}_A\| \|e_b\|,\\
V_3(\|\mathcal{E}_A\|, \|e_b\|) &= (k_P -(2B_\xi + B_b)- \epsilon k_I) \|\mathcal{E}_A\|^2 \\
&\quad+ \epsilon \|e_b\|^2  -\epsilon (k_P + B_b+2B_\xi)\|E_A\| \|e_b\|.
\end{align*}
Hence, there are numbers $\alpha>0$ and $\beta >0$ such that  \eqref{V213} holds. 
Let 
\[
V(\mathcal{E}_A, e_b) = \frac{1}{2} \|\mathcal{E}_A\|^2 + \frac{1}{2k_I} \|e_b\|^2 + \epsilon \langle \mathcal{E}_A,  e_b \rangle,
\]
which satisfies \eqref{V1:V:V2}
for all $(\mathcal{E}_A, e_b) \in \mathbb R^{n\times n} \times \Lalge$. It is easy to show that along any trajectory of the composite system consisting of the rigid body \eqref{rigid:body} and the observer \eqref{observer:Chang:three},
\[
\dot V  \leq -V_3 \leq -\beta V_2 \leq -\beta V.
\]
As in the proof of Theorem \ref{main:theorem:observer}, it is east to show that there are numbers $\tilde C>0$ and $ a>0$ such that 
\begin{equation}\label{append:three:1}
\| {\mathcal E}_A(t)\| + \|e_b(t)\| \leq \tilde C( \| {\mathcal E}_A(0)\| + \|e_b(0)\| )e^{-at}
\end{equation}
for all $(\bar A(0), \bar b(0)) \in \mathbb R^{n\times n} \times \Lalge$ and all $t\geq 0$.  Since ${\mathcal E}_A = I - A^{-1}\bar A = A^{-1}E_A$ or $E_A = A {\mathcal E}_A = \meas g {\mathcal E}_A$,  we have
\begin{equation}\label{append:three:2}
 \sigma_{\rm min}(\meas ) L_g \|{\mathcal E}_A\| \leq \|E_A\| \leq \sigma_{\rm max}(\meas ) U_g \|{\mathcal E}_A\|.
\end{equation}
It follows from \eqref{append:three:1} and \eqref{append:three:2} that there is a number $C>0$ such that 
\[
\| E_A(t)\| + \|e_b(t)\| \leq  C( \| E_A(0)\| + \|e_b(0)\| )e^{-at}
\]
for all $(\bar A(0), \bar b(0)) \in \mathbb R^{n\times n} \times \Lalge$ and all $t\geq 0$. 
\end{proof}


\end{document}